\documentclass[12pt]{article}

\usepackage{amssymb,amsmath,amscd}

\newtheorem{df}{Definition}[section]
\newtheorem{thm}{Theorem}[section]

\newtheorem{rem}{Remark}[section]

\def\GA{\mathcal{G}_{V}(A)}
\def\GB{\mathcal{G}_{V}(B)}

\def\G{\mathcal{G}_{V}}
\def\GC{G_{\mathbb{C}}}

\begin{document}


\vspace{3cm}
\begin{center}
{\Large \bf On real forms of complex Lie superalgebras and complex algebraic supergroups}\\
[50pt]{\small
{\bf F. Pellegrini}\,\footnote{e-mail:pelleg@iml.univ-mrs.fr}
\\ ~~\\Institut de math\'ematiques de Luminy,
 \\ 163, Avenue de Luminy, 13288 Marseille, France}


\vspace{1cm}
\begin{abstract}
The paper concerns two versions of the notion of real forms of Lie superalgebras. One is the standard approach, where a real form of a complex Lie superalgebra is a real Lie superalgebra such that its complexification is the original complex Lie superalgebra. The second is related to considering $A$-points of a Lie superalgebra over a commutative complex superalgebra $A$  equipped with superconjugation. It is not difficult to see that the first real form can be obtained as the set of fixed points of an antilinear involutive automorphism and the second is related to an automorphism $\phi$ such that $\phi^{2}$ is identity on even part and negative identity on the odd part. The generalized notion of the real form is subsequently introduced also for complex algebraic supergroups. 
\end{abstract}
\end{center}

\newpage


\section{Introduction}
  
There are two equivalent points of view to define a complex Lie superalgebra. The standard one is about complex Lie superalgebra as a supervector space with a superbracket. The other defines Lie superalgebra as a representable functor. The equivalence of the two definitions is a corollary of the so called "even rules" \cite{Deligne}. In studying any particular aspect of the theory of Lie superalgebras, we have a choice to work either with the standard (supervector) Lie superalgebra or with the functorial definition. We shall argue in this paper, that the functorial language is particularly well adapted to the problem of definition and classification of real forms of simple Lie superalgebras.
\vskip1pc 
\noindent The real forms were defined and classified by Serganova \cite{Serganova} in the standard framework of complex vector superspaces with a somewhat puzzling conclusion that the simple complex Lie superalgebras have no compact real forms. It is the subject of this paper to argue that the functorial point of view avoids this puzzle by defining two kinds of real forms which we call standard and graded one. In the standard approach the real form is defined as a real Lie superalgebra such that its complexification is the original complex Lie superalgebra. It can be seen easily that every standard real form is naturally associated to an antilinear involutive automorphism of the complex Lie superalgebra \cite{Serganova}. We shall see that the notion of graded real form is characterised by antilinear automorphisms which are not involutive but rather graded involutive. It turns out, remarkably, that Serganova has classified in her paper also the graded involutive automorphisms although she did not interpret them as real forms. Using the Serganova results, we show that our more general definition of real forms solves the puzzle, in the sense that every simple complex Lie superalgebra turns out to possess precisely one compact (graded) real form.
\vskip1pc
\noindent The functorial point of view turns out to be useful also for defining the real forms of complex algebraic supergroups. We do it for the subsupergroups of the $GL_{m \vert n}$. We show also that a real form of a supergroup induces  a real form on its corresponding Lie superalgebra. This rises a question what are the lifts of the Serganova automorphisms to the corresponding supergroups. 
We answer it completely for the series $SL_{m \vert n}$ and $OSp_{m \vert 2n}$.


\section{Graded and standard real forms of complex Lie superalgebras}

\noindent In this work a "superalgebra" will stand always for a commutative and associative superalgebra over $\mathbb{C}$. The even and odd part of a superalgebra $A$ are respectively denoted $A_{0}$ and $A_{1}$. Elements of $A_{0}, A_{1}$ are said homogeneous and the parity of an homogeneous element $x$ is $\vert x \vert=0$ for $x \in A_{0}$ or $\vert x \vert=1$ for $x \in A_{1}$. A superalgebra is said commutative if $xy=(-1)^{\vert x \vert \vert y \vert}yx$ for homogeneous elements $x,y$. $\mathcal{I}_{A}^{odd}$ is the ideal in $A$ generated by the odd part. We say that a superalgebra $A$ is reduced if $A/\mathcal{I}_{A}^{odd}$ have no nilpotent elements. 
We denote by $\mathcal{A}$ and $\mathcal{S}$ respectively the category of reduced, finetely generated complex superalgebras (such superalgebra is also called affine superalgebra) and the category of sets. We restrict our category $A$ to the reduced and finitely generated superalgebra because in this article we are only concerned with the Lie superalgebra which come from affine algebraic supergoup.
Let $\mathcal{O}: \mathcal{A} \longrightarrow \mathcal{S}$ be the functor which sends each superalgebra $A$ on its even part i.e. $\mathcal{O}(A)=A_{0}$. 
Now, we present the definition of complex Lie superalgebra in the functorial setting as it is given in \cite{Fioresi}:  


\begin{df}

A Lie superalgebra is a representable functor 
\begin{eqnarray}
\G :\mathcal{A} &  \longrightarrow & \mathcal{S} \\
A  & \longmapsto & \G (A)
\end{eqnarray}
such that:
\vskip1pc
\noindent 
i)for each superalgebra $A$ we have an $A_{0}$-module structure on $\G(A)$ which is functorial in $A$ in other words there is a natural transformation $\mathcal{O} \times \G \rightarrow  \G$, 
\vskip1pc
\noindent 
ii)for each superalgebra $A$, we have a Lie-bracket $[\, ,\, ]_{A}$ on $\G(A)$ which is $A_{0}$-linear and functorial in $A$ in other words there is a natural transformation $[\, , \, ] : \G \times \G \rightarrow \G$ which is $\mathcal{O}$-linear and satisfies commutative diagram corresponding to the antisymmetric property and Jacobi identity.

\end{df}


\vskip1pc
\noindent
The statement that the Lie superalgebra is representable means that it exists a supervector space $V$ (which is unique) such that $\G(A)=(A \otimes_{\tiny{\mathbb{C}}} V)_{0}$. Following the even rules principle, $V$ inherits a superbracket from the Lie algebra structure of the $\G(A)$ for all $A$. For clarity,  we explain the even rules principle (for a proof \cite{Deligne} Theorem 1.7.1 pg 56 and corollary 1.7.3 pg 57 and \cite{Varadarajan}) in the extent needed for this paper. If $f$ is a bilinear map from $V \times V$ to $V$ then  for $b_{1}, b_{2}\in B, \,v_{1}, v_{2} \in V$, which are such that $\vert b_{i} \vert =\vert w_{i} \vert$, and $N$ is the numbers of  odd elements among the $b_{i}$, we define a $B_{0}$-bilinear map $f_{B}:\GB \times \GB \rightarrow \GB$ by 
\begin{eqnarray}
f_{B}(b_{1} \otimes v_{1}, b_{2} \otimes v_{2})=(-1)^{\frac{N(N-1)}{2}}b_{1}b_{2} \otimes f(v_{1}, v_{2}).
\end{eqnarray} 
It is not difficult to observe that the collection of maps $f_{B}$ is functorial in $B$. This means that the following diagram is commutative for any complex superalgebras $B,C$
\begin{equation}
\begin{CD}
\GB \times \GB @>f_{B}>> \GB\\
@V (\G(f), \G(f)) VV     @VV \G(f) V\\
\GC \times \GC @>f_{C}>> \GC
\end{CD}
\end{equation} 
where the morphism $\G(f)$ associated to $f:B \rightarrow C$ is given by 
$\G(f)(b \otimes v)=f(b)\otimes v$ for $b \in B,\, v\in V$.
Now the principle of even rules claims that any functorial collection $f_{B}$ of $B_{0}$-linear maps 
\begin{eqnarray}
f_{B}:\GB \times \GB \rightarrow \GB
\end{eqnarray} 
comes from a unique $\mathbb{C}$-linear map $f: V \times V  \rightarrow V.$ i.e. the maps $f_{B}$ and $f$ are linked by the formula $(1)$. 
Thus a Lie superbracket $[\,,\,]$ comes from the functorial collection of Lie bracket $[\,,\,]_{A}$ and , vice versa, a Lie superbracket gives rise to a system of functorial brackets. The relation between $[\, , \,]$ and $[\, , \,]_{A}$ is 
\begin{eqnarray}
[a \otimes v, b \otimes w]_{A}=(-1)^{\vert b \vert \vert v \vert}ab \otimes [v,w],
\end{eqnarray}   
where $a,b \in A ,\, v,w\in V$ which fullfill $\vert a \vert =\vert v \vert$ and $\vert b \vert = \vert w \vert$.
\vskip1pc
\noindent
Now, we want to introduce the notion of real structure of a complex Lie superalgebra. For this, the category $\mathcal{A}$ from the definition $2.1$ must be constructed as the category of superalgebras with conjugation. In particular, this means that the morphism in $\mathcal{A}$ respects the conjugation (i.e. the morphism commute with the conjugation). In distinction to the non-super case there is two kinds of conjugation on the commutative associative superalgebra $\cite{Deligne}$:


\begin{df}
Let $A$ be a complex superalgebra. A map $a \rightarrow \bar{a}$ with $a, \bar{a}\in A$ is called a \underline{standard} conjugation if it holds:
\begin{eqnarray}
\quad \overline {\lambda a}=\bar \lambda \bar a,
\quad \quad \overline{a b}= \bar a \bar b, \quad  \quad \bar{ \bar a} = a
,  \quad \lambda \in \mathbb{C}, \quad a,b \in A.
\end{eqnarray}
\end{df}  

\begin{df}
Let $A$ be a complex superalgebra. A map $a \rightarrow \tilde{a}$ with $a, \tilde{a}\in A$ is called a \underline{graded} conjugation if it holds:
\begin{eqnarray}
\quad \widetilde{\lambda a} = \bar{\lambda} \tilde{a},
\quad \quad \widetilde{ab}= \tilde{a} \tilde{b}, 
\quad  \quad \tilde{ \tilde{a}} = (-1)^{\vert a \vert} a,  \quad \lambda \in \mathbb{C}, \quad a,b \in A.
\end{eqnarray}
\end{df}



\begin{rem}
We denote the conjugated element of $a \in A$ as $\hat{a}$ when the distinction between the graded and the standard conjugation does not matter. We say that an element $a \in A$ is real if it fullfills $a=\hat{a}$. The superalgebra $A^{real}$ is the superalgebra consisting of the real elements of $A$.
\end{rem}

 
\vskip1pc
\noindent
Before to give the definition of the real structure of a Lie superalgebra in the functorial setting we have to give some conventions. In the case of a supervectorspace $V=V_{0}+V_{1}$ a map $f:V\rightarrow V$ is said even when $f$ sends respectively the elements of $V_{0}$ ($V_{1}$) in $V_{0}$ ($V_{1}$) i.e. it respects the parity of $V$. The $A_{0}$-module $(A \otimes V)_{0}$ have also a gradation, it is given by the following decomposition 
\begin{eqnarray}
(A \otimes V)_{0} = A_{0} \otimes V_{0}+A_{1} \otimes V_{1}.
\end{eqnarray}
The elements of $A_{0} \otimes V_{0}$ ($A_{1} \otimes V_{1}$) are said evens (odds). A map f on $(A \otimes V)_{0}$ is  even when it preserves the previous gradation.


\begin{df}
A real structure of a Lie superalgebra is a natural transformation $\Phi : \mathcal{G} \rightarrow \mathcal{G}$ such that for $a, b \in A_{0},\, x,y \in \mathcal{G}(A)$:
\begin{eqnarray}
\Phi_{A}(a x + b y) & = & \hat{a} \Phi_{A}(x) + \hat{b} \Phi_{A}(y),\\
\Phi_{A}(\Phi_{A}(x)) & = & x, \\
\Phi_{A}([x,y]_{A}) & = & [\Phi_{A}(x),\Phi_{A}(y)]_{A},
\end{eqnarray}
Furthermore we demand that $\Phi_{A}$ is an even map.
The real structure is called \underline{standard} (\underline{graded}) when the conjugation is the standard (graded) one.  
\end{df}



\begin{thm}
1) A standard real structure $(\G,\Phi)$ comes from a unique map $\phi:V \rightarrow V$ such that
\begin{eqnarray}
\phi(\lambda x + \mu y)= \bar{\lambda} \phi(x) + \bar{\mu} \phi(y),\, \phi^{2}(x)=x,\, \phi([x,y])=[\phi(x),\phi(y)].
\end{eqnarray}
with $\lambda, \mu \in \mathbb{C},\, x,y \in V$.
\vskip1pc
\noindent
A graded real structure $(\G,\Phi)$ comes from a unique map $\phi:V \rightarrow V$ such that
\begin{eqnarray}
\phi(\lambda x + \mu y)= \bar{\lambda} \phi(x) + \bar{\mu} \phi(y),\, \phi^{2}(x)=(-1)^{\vert x \vert}x,\, \phi([x,y])=[\phi(x),\phi(y)].
\end{eqnarray} 
with $\lambda, \mu \in \mathbb{C},\, x,y \in V$.
\vskip1pc
\noindent
2) Reciprocally, a map $\phi$ on $V$ which fullfills the conditions $(11)$ (resp. $(12)$) gives rise to a standard (resp. graded) real form $(\G,\Phi)$.
\end{thm}


\vskip1pc
\noindent
\textbf{Proof:}
We give the proof of this theorem for the graded real form. The proof in the case of the standard real form is very similar.


\vskip1pc
\noindent 1) Let $(\G,\Phi)$ be a graded real form. First we want to show that $\Phi$ gives rise to a map $\phi: V \rightarrow V$, which means that $\Phi(a \otimes v)=\tilde{a} \otimes \phi(v)$ with $a \in A,\, v \in V$ and $\vert a \vert = \vert v \vert$. The crucial property used is the fact that $\Phi$ is a natural transformation. Thus for $A,B$ two superalgebras and the morphism of superalgebra $p: A \rightarrow B$ we have the following commutative diagram
\begin{equation}
\begin{CD}
\G(A) @>\Phi_{A}>> \G(A)\\
@V \G(p) VV     @VV \G(p) V\\
\G(B) @>\Phi_{B}>> \G(B).
\end{CD}
\end{equation}
We say that this diagram is associated to the triple $(A,B,p)$. The morphism $p$ fixes clearly the orientation of the diagram. First we want to extract a map $\phi: V \rightarrow V$ from $\Phi_{A}$ for $A=\mathbb{C}[\theta,\tilde{\theta}]$\, 
\footnote{$A$ is an example of Grassmann algebra with a graded conjugation. A Grassmann algebra is a superalgebra $\mathbb{C}[\theta_{1},...,\theta_{n}]$ over $\mathbb{C}$ generated by $n$ odd elements $\theta_{i}$ , which anticommute i.e. $\theta_{i}\theta_{j}+\theta_{j}\theta_{i}=0$ for all $i$, in particular we have $\theta_{i}^{2}=0$. The notation $\mathbb{C}[\theta,\tilde{\theta}]$ clearly indicates how the graded conjugation acts on the generators.}
. For this we have to evaluate $\Phi_{A}$ only on $1 \otimes v$ and $\theta \otimes w$ with $\vert v \vert=0,\, \vert w \vert=1,\, \lambda \in \mathbb{C}$. Note that it is not necessary to compute also $\Phi_{A}(\tilde{\theta}\otimes w))$ because the graded conjugation is a superalgebra morphism and $\Phi$ is a natural transformation. In order to evaluate $\Phi_{A}(1 \otimes v)$ we use the commutative diagram associated to the triplet $(\mathbb{C},A,p)$ where $p$ is the morphism which injects the complex number. We find that $\Phi_{A}(1 \otimes v)= 1 \otimes \Phi_{\mathbb{C}}(v)$. Now let $\pi: A \rightarrow \mathbb{C}$ be the canonical projection $\pi(\theta)=\pi(\tilde{\theta})=0$. From the commutative diagram associated to $(A,\mathbb{C},\,\pi)$ and the fact that $\Phi_{A}$ is an even map we deduce that it exists $x,y \in V_{1}$ such that
\begin{eqnarray}
\Phi_{A}(\theta \otimes w)= \theta \otimes x + \tilde{\theta} \otimes y.
\end{eqnarray}
Moreover, in defining the morphism of superalgebra $q: A \rightarrow D$ with $A=\mathbb{C}[\theta,\tilde{\theta}],\, D=\mathbb{C}[\theta,\tilde{\theta},\eta,\tilde{\eta}]$ (such that $q(\theta)=\theta,q(\tilde{\theta})=\tilde{\theta}$), we obtain from the diagram $(A,D,q)$: 
\begin{eqnarray}
\Phi_{D}(\theta \otimes w)= \theta \otimes x +\tilde{\theta} \otimes y.
\end{eqnarray}
But by multiplying this last equation by $\eta \tilde{\theta}$ we find $\eta \tilde{\theta}\theta \otimes x =0$ therefore $x=0$. Thus 
\begin{eqnarray}
\Phi_{D}(\theta \otimes w)= \tilde{\theta} \otimes y.
\end{eqnarray}
Then from the commutative diagram associated to the triplet $(D,A,\Pi)$, where $\Pi$ is defined by $\Pi(\eta)=\Pi(\tilde{\eta})=0$, we deduce 
\begin{eqnarray}
\Phi_{A}(\theta \otimes w)=\tilde{\theta} \otimes y 
\end{eqnarray}
Now we define $\phi$ by $\phi(w)=y$ and $\phi(v)=\Phi_{\mathbb{C}}(v)$ for $v \in V_{0}, \, w \in V_{1}$. We have to prove that this map $\phi$ gives rise to $\Phi_{C}$ for any superalgebra $C$. In others words, we have to show that $\Phi_{C}(c \otimes v )= \tilde{c} \otimes \phi(v)$ for $\vert c \vert = \vert v \vert$. For $c \in C_{0}$ the commutative diagram associated to $( \mathbb{C}[\theta,\tilde{\theta}],C, m(1)=1 )$ imply $\Phi_{C}(1 \otimes v)=1 \otimes \phi(v)$ and thus by $C_{0}$-antilinearity we have $\Phi_{C}(c \otimes v)=\tilde{c}\otimes \phi(v)$. For $c \in C_{1}$ we define the following morphism of superalgebra $n(\theta)=c, \, n(\tilde{\theta})=\tilde{c}$ from $\mathbb{C}[\theta,\tilde{\theta}]$ to $C$ and we obtain the equality desired by the commutativity of the diagram $(\mathbb{C}[\theta,\tilde{\theta}],C,n)$. The final step is to show that $\phi$ fullfills the properties $(11)$. The antilinearity is clear from the $C_{0}$-antilinearity of $\Phi_{C}$ (in particular the $\mathbb{C}$-antilinearity) and the linearity of the tensor product. From the involutivity of $\Phi_{C}$ and the definition of $\phi$ ($\Phi(c \otimes v)=\tilde{c}\otimes \phi(v)$) we find $\tilde{\tilde{c}}\otimes \phi^{2}(v)=c \otimes v$. Therefore as $\vert c \vert = \vert v \vert$ we have $\phi^{2}(v)=(-1)^{\vert v \vert}v$. The fact that $\phi$ is a morphism of Lie superalgebra comes from the property of morphism of Lie algebra of $\Phi_{C}$ and the equation $(3)$. Thus we have ahieved the demonstration of $1)$.


\vskip1pc
\noindent 2) Let $(V,\phi)$ be a graded real form. For each superalgebra $A$ we define the collection of maps $\Phi_{A}: \GA \rightarrow \GA$ by the formula
\begin{eqnarray}
\Phi_{A}(a \otimes v)= \tilde{a} \otimes \phi(v)
\end{eqnarray} 
for all $a \in A$ and $v \in V$. Thus, from this formula, we find that
\begin{eqnarray}
\G(f)(\Phi_{A}(a \otimes v))=\Phi_{B}(\G(f)(a \otimes v)).
\end{eqnarray}
This equation comes from the definition of $\G(f)$ and the equality $f(\tilde{a})=\widetilde{f(a)}$. The equality $(19)$ means that the collection of maps $\Phi_{A}$ is functorial in $A$. The property of antilinearity in $A_{0}$ comes from the definition of $\Phi_{A}$ in $(18)$, the linearity of the tensor product and the  structure of $A_{0}$-module of $\GA$.
Then we have from $(18)$ that for $a \in A, \, v \in V$ with $\vert a \vert =\vert v \vert$
\begin{eqnarray}
\Phi_{A}(\Phi_{A}(a \otimes v))=\Phi_{A}(\tilde{a} \otimes \phi(v))=\tilde{\tilde{a}}\otimes \phi(\phi(v))=(-1)^{\vert a \vert +\vert v \vert} a \otimes v = a \otimes v.
\end{eqnarray}
Thus $\Phi_{A}$ is involutive. From the fact that $\phi$ is a morphism of Lie superalgebra  and the equation $(3)$ we deduce easily that $\Phi_{A}$ is a morphism of the Lie algebra $\G(A)$. Hence $\Phi$ is a graded real structure.
$\qquad\qquad\qquad\qquad\qquad\qquad \blacksquare$


\vskip1pc
\noindent
In the \underline{standard} setting we can associate to each \underline{real structure} $\Phi$ the corresponding \underline{real form} which is by definition the real superalgebra $V^{\phi}=\{v \in V,\, \phi(v)=v \}$ obtained as the fixed point set of the automorphism $\phi$. However the corresponding notion of the graded real form as the fixed point set of $\phi$ is more subtle because it appears that $V^{\phi}=\{v \in V,\, \phi(v)=v \}$ is a trivial Lie superalgebra (i.e. the odd elements were killed by the requirement $\phi(v)=v$). It turns out that the correct point of view which resolve this trouble is again functorial since it treats the notion of the real form on the same footing for both standard and graded real structure. Thus we give the following definition:


\begin{df}
Let $\Phi$ be a real structure of $\G$. The real form associated to $\Phi$ is the functor $\G^{\Phi}$ 
\begin{eqnarray*}
\G^{\Phi}: \mathcal{A} &  \longrightarrow & \mathcal{S} \\
A  & \longmapsto & \G^{\Phi}(A)
\end{eqnarray*}
where $\G^{\Phi}(A)=\{ w \in \GA, \Phi_{A}(w)=w \}$ is a $A_{0}^{real}$-Lie algebra and to a morphism $f :A \rightarrow B$ of superalgebras it associates $\G^{\Phi}(f)=\G(f)_{\G^{\Phi}(A)}$ (this means that $\G^{\Phi}(f)$ is the restriction of $\G(f)$ to the set of fixed points of $\Phi_{A}$). The real form is said standard (graded) when the real structure is standard (graded).
\end{df}
 
 

\begin{rem}
The definition 2.5 is consistent because it is easy to demonstrate that $\G^{\Phi}$ is a functor. Moreover, $\G^{\Phi}(f)$ sends the elements of $\G^{\Phi}(A)$ into $\G^{\Phi}(B)$ because $\Phi$ is a natural transformation.
\end{rem}


The fact that we can extract a real Lie superalgebra $V^{\phi}$ from a standard real structure while it is impossible to do it for a graded real structure, have some consequences on the representability of the real form. In fact we have the following theorem.


\begin{thm}
1) If $\Phi$ is a \underline{standard} real structure then the functor $\G^{\Phi}$ is represented by $V^{\phi}$ i.e.  $\G^{\Phi}(A)=(A^{real}\otimes V^{\phi})_{0}$.\\
\noindent
2) If $\Phi$ is a \underline{graded} real struture then the functor $\G^{\Phi}$ is not representable.
\end{thm}

 
\vskip1pc
\noindent
\textbf{Proof:}
1) We begin with the proof of the first implication. Thus, $A$ is equipped with a standard conjugation and $\Phi$ is a standard real structure.
Let $a \otimes v \in \G(A)$ be such that $\Phi_{A}(a \otimes v)= a \otimes v$.
We know from the theorem 2.1 that a standard real structure $\Phi$ comes from a unique map $\phi$ which fullfills the properties $(9)$. We remark that every element $v$ of $V$ can be uniquely decomposed as $v=v_{1}+iv_{2}$ with $\phi(v_{1})=v_{1}$ and $\phi(v_{2})=v_{2}$. Similarly, every elements $a$ of $A$ have the unique decomposition $a=a_{1}+ia_{2}$ with $\bar{a}_{1}=a_{1}$ and $\bar{a}_{2}=a_{2}$. Thus we have the following equalities:
\begin{eqnarray*}
a \otimes v  =  \frac{a \otimes v +\Phi_{A}(a\otimes v)}{2}
 =  a_{1}\otimes v_{1} - a_{2} \otimes v_{2}
\end{eqnarray*} 
From this last equality, it is clear that the fixed points of $\Phi_{A}$ are elements of $(A^{real}\otimes V^{\phi})_{0}$.
\vskip1pc
\noindent
2) Now $A$ and $V$ be equipped, respectively, with a graded conjugation $a \rightarrow \tilde{a}$ and a map $\phi$ which comes from a graded real structure $\Phi$ on $\G$. It is clear that odd element of $A$ cannot be real for a graded conjugation and  odd vector of $V$ cannot be a fix points of $\phi$. So, there are some fixed points of $\Phi_{A}$ in $\G(A)$ which are not elements of  $(A^{real}\otimes V^{\phi})_{0}$, for instance, the following one $a \otimes v +\tilde{a} \otimes \phi(v)$ with $a \in A_{1}$ and $v \in V_{1}$. $\qquad\qquad\qquad\qquad\qquad\qquad\qquad\blacksquare$


\vskip1pc
\noindent
Denote as $(\G^{\Phi})_{0}$ the restriction of the functor $\G^{\Phi}$ to its even part with respect to the gradation  define by the decomposition $(6)$, i.e. ${(\G^{\Phi}})_{0}(A)=(\G^{\Phi}(A))_{0}$. It turns out that this restricted functor is representable for both standard and graded real forms (the proof is similar to the demonstration of the theorem 2.2). The representative of ${(\G^{\Phi}})_{0}$ is an ordinary Lie algebra $(V^{\phi})_{0}$. It makes therefore sense to ask whether $(V^{\Phi})_{0}$ is compact. Remarkably, it follows from the Serganova classification (see table 3 and table 6 in \cite{Serganova}) that there is exactly one (graded) automorphism  $\phi$ for each simple complex Lie superalgebra such that $(V^{\phi})_{0}$ is always compact, which shows the utility of our interpretation of graded involutive automorphisms as real forms.


\section{Graded and standard real forms of complex algebraic supergroup}

The goal of this section is to introduce the notions of real structure and real form for subsupergroup of $GL_{m \vert n} $. 
\vskip1pc
\noindent
In order to define the supergroup $GL_{m \vert n}$, we have to introduce some notations. We denote by $A^{m \vert n}$ the free $A$-supermodule generated by $m$ even $\mathbf{e}_{1},...,\mathbf{e}_{m}$ and $n$ odd generators $\mathbf{e}_{m+1},...,\mathbf{e}_{m+n}$ such that $\mathbf{a} \in A^{m \vert n}$ is of the form $\mathbf{a}=a_{1}\mathbf{e}_{1}+...+a_{m+n}\mathbf{e}_{m+n}$. An even morphism $T: A^{m \vert n}\rightarrow A^{m \vert n}$ can be represented by a supermatrix $\mathcal{T}$ of size $(m+n) \times (m+n)$ 
\begin{eqnarray}
\mathcal{T}=\left(
\begin{array}{cc}
P & Q \\
R & S 
\end{array}
\right) 
\end{eqnarray}
where the matrices $P,S$ have even entries from $A$ and are respectively of size $m \times m$, $n \times n$; the matrices $Q,R$ have odd entries and are respectively of size $n \times m$, $m \times n$. The endomorphisms and the automorphisms of $A^{m \vert n}$ are denoted respectively $gl_{m \vert n}(A)$ and $GL_{m \vert n}(A)$. The supermatrices $\mathcal{T}$ of $GL_{m \vert n}(A)$ are such that the Berezinian or superdeterminant (see \cite{Berezin})
\begin{eqnarray}
sdet(\mathcal{T})=det(P-QS^{-1}R)det(S^{-1}), 
\end{eqnarray} 
is invertible in $A$. A necessary and sufficient condition for invertibility of $sdet(\mathcal{T})$ is the invertibility of $P$ and $S$. The functor of Linear affine algebraic supergroup is
\begin{eqnarray}
GL_{m \vert n}:\mathcal{A} &  \longrightarrow & \mathcal{S}  \\
A  & \longmapsto & GL_{m \vert n}(A). \nonumber
\end{eqnarray}


\begin{df}
A complex linear affine algebraic subsupergroup (in the sequel we call it just complex algebraic supergroup) is a functor $G_{\mathbb{C}}$
\begin{eqnarray}
G_{\mathbb{C}}:\mathcal{A} &  \longrightarrow & \mathcal{S} \\
A  & \longmapsto & G_{\mathbb{C}}(A) \nonumber
\end{eqnarray}
where $G_{\mathbb{C}}(A)=\{x \in GL_{m\vert n}(A), \; \mathcal{P}_{l}(x)=0, \;  \forall l=1...k \}$ and the polynomials $\mathcal{P}_{l}$ are such that $G_{\mathbb{C}}(A)$ is a group.
\end{df}


\begin{rem}
Two examples are the $SL_{m\vert n},\, OSp_{m\vert 2n}$ series of supergroup where the polynoms $\mathcal{P}_{l}$ are defined in the equation $(40), (41)$ respectively.
\end{rem}

\noindent Let $A$ be a superalgebra and $\epsilon$ a formal inderterminate. Let $A(\epsilon)$ be the superalgebra of dual numbers defined by $A(\epsilon)=(A\oplus \epsilon A)/(\epsilon^{2})$. There are three useful morphisms  $i: A \rightarrow A(\epsilon)$ defined by $i(x)=x+ \epsilon \, 0$; $p:A(\epsilon) \rightarrow A$ defined by $p(x+\epsilon \, y)= x$ and $v_{a}:A(\epsilon) \rightarrow A(\epsilon)$ defined by $v_{a}(x+\epsilon \, y)=
x+ \epsilon \, a y$ for $a \in A_{0},\, x,y \in A $.   


\begin{df}
A real structure of a complex algebraic supergroup $\GC$ is a natural transformation $\Sigma: \GC \rightarrow \GC$ which fullfils
\begin{eqnarray}
\; \Sigma_{A}(xy) &=& \Sigma_{A}(x)\Sigma_{A}(y), \\
\; {\Sigma_{A}}^{2}(x) &=& x \\
\Sigma_{A(\epsilon)}(\GC(v_{a})(z)) &=& \GC(v_{\hat{a}})(\Sigma_{A(\epsilon)}(z))
\end{eqnarray}
for all $x,y \in \GC(A),\; z \in Ker(\GC(p))$. The real structure is said standard (graded) when the superalgebra A is equipped with a standard (graded) conjugation.
\end{df}


\begin{rem}
The requirement $(29)$ means that the map induced by $\Sigma$ on the Lie superalgebra of $\GC$ is antilinear.
\end{rem}


\begin{df}
Let $\Sigma$ be a real \underline{structure} of $\GC$. The real \underline{form} associated to $\Sigma$ is the functor 
\begin{eqnarray}
\GC^{\Sigma}: \mathcal{A} &  \longrightarrow & \mathcal{S} \\
A  & \longmapsto & \GC^{\Sigma}(A)
\end{eqnarray}
where $\GC^{\Sigma}(A)= \{ x \in \GC(A), \, \Sigma_{A}(x)=x \}$ is a group and to a morphism $f: A \rightarrow B$ of superalgebras it associates $\GC^{\Sigma}(f)=\GC(f)_{\GC^{\Sigma}}$ (this means that $\GC^{\Sigma}(f)$ is the restrication of $\GC(f)$ to the set of fixed points of $\Sigma_{A}$). The real form is said standard (graded) when the real structure is standard (graded).
\end{df}


\noindent Now we define the Lie superalgebra associated to a supergroup (see \cite{Fioresi}).

\begin{df}
The Lie superalgebra of the algebraic supergroup $\GC$ is the functor
\begin{eqnarray}
Lie(\GC): \mathcal{A} &  \longrightarrow & \mathcal{S} \\
A  & \longmapsto & Lie(\GC)(A)= Ker(\GC(p))
\end{eqnarray}
where $\GC(p): \GC(A(\epsilon)) \rightarrow \GC(A)$ is the morphism associated to $p:  A(\epsilon) \rightarrow A$. 
\end{df}

\begin{rem}
As explained in $\cite{Fioresi}$ $Ker(\GC(p))$ have a structure of $A_{0}$-module with a Lie-bracket.
\end{rem}

\noindent From the definition 3.1 we deduce that the set $Ker(\GC(p))$ is the set of even supermatrices $N$ such that $\mathcal{P}_{l}(1+ \epsilon N)=0$ for all $l=1...k$. From \cite{Fioresi} we deduce that the Lie-bracket on $Ker(\GC(p))$ is simply the commutator $[M,N]=MN-NM$. Moreover $Ker(\GC(p))$ is an abelian group with respect to the addition of supermatrices and the action $A_{0}$ on $Ker(\GC(p))$ corresponds to the multiplication of the entries of   the supermatrices by an element of $A_{0}$. Thus $Ker(\GC(p))$ is also a $A_{0}$-module. 
\vskip1pc
\noindent Now we prove that the real structure of a supergroup gives rise to a real structure of its corresponding Lie superalgebra.


\begin{thm}
Let $\GC$ be a algebraic supergroup and $\Sigma$
be a real structure of $\GC$ then for every superalgebra $A$ we have 
\begin{eqnarray}
\Sigma_{A(\epsilon)}(Id+\epsilon M)=Id+\epsilon \Phi_{A}(M),\; \forall M \in \GC(A).
\end{eqnarray}
The collection of maps $\Phi_{A}$ define a real structure on $Lie(\GC)$.
\end{thm}


\vskip1pc
\noindent
\textbf{Proof:}
Let $\GC(p): \GC(A(\epsilon))\rightarrow \GC(A)$ be the morphism induced by the morphism $p: A(\epsilon)\rightarrow A$ such that $p(a+\epsilon b)=a,\; \forall a,b \in A$. $\Sigma_{A(\epsilon)}$ is a morphism from $\GC(A(\epsilon))$ to $\GC(A)$ hence, by functoriality, we have $\GC(p)\Sigma_{A(\epsilon)}(Id+ \epsilon M)=\Sigma_{A}(\GC(p)(Id+ \epsilon M))=\Sigma_{A}(Id)=Id$. We deduce that $\Sigma_{A(\epsilon)}(Id+ \epsilon M)$ is an element of $Ker(\GC(p))$ therefore we have defined a map $\Phi_{A}$ such that $\Sigma_{A(\epsilon)}(Id+\epsilon M)=Id + \Phi_{A}(M)$. In order to show that $\Phi$ is a real structure of $Lie(\GC)$, we have to prove that each maps $\Phi_{A}$ is antilinear, involutive and is a morphism of Lie algebra. We conclude that $\Sigma_{A(\epsilon)}(Id + \epsilon (M+N))= Id +\epsilon \Phi_{A}(M+N)$ and $\Sigma_{A(\epsilon)}(Id + \epsilon a M) = Id +\Phi_{A}(aM)$ with $a\in A_{0}$ and $M,N \in \GC(A)$. Furthermore we have
\begin{eqnarray}
\Sigma_{A(\epsilon)}(Id + \epsilon (M+N)) & = &\Sigma_{A(\epsilon)}((Id+\epsilon M)(Id+\epsilon N))\nonumber\\
& = & \Sigma_{A(\epsilon)}(Id+\epsilon M)  \Sigma_{A(\epsilon)}(Id+\epsilon N)\nonumber\\
& = & (Id + \epsilon \Phi_{A}(M))(Id + \epsilon \Phi_{A}(N))\nonumber\\
& = & Id + \epsilon (\Phi_{A}(M)+ \Phi_{A}(N))
\end{eqnarray}
and
\begin{eqnarray}
\Sigma_{A(\epsilon)}(Id + \epsilon a M) & = & \Sigma_{A(\epsilon)}(\GC(v_{a})(Id + \epsilon M)))\nonumber\\
& = & \GC(v_{\hat{a}})(\Sigma_{A(\epsilon)}(Id + \epsilon  M))\nonumber\\
& = & \GC(v_{\hat{a}})(Id + \epsilon \Phi_{A}(M))\nonumber\\
& = & Id+\epsilon \hat{a} \Phi_{A}(M).
\end{eqnarray}
The equations $(35)$ and $(36)$ establish the antilinearity of $\Phi_{A}$. The involutivity of $\Phi_{A}$ is deduced easily from the involutivity of $\Sigma_{A(\epsilon)}$.\\
It remains to show that $\Phi_{A}$ is a morphism of Lie-algebra. Let $A(\epsilon, \eta)$ be the superalgebra of polynomials of the indertermine $\epsilon$ and $\eta$ with coefficients in $A$ and such that $\epsilon^{2}=0$, $\eta^{2}=0$, $\epsilon \eta -\eta \epsilon=0$. We have two morphisms $p_{\epsilon}:A(\epsilon,\eta)\rightarrow A(\eta)$ and $p_{\eta}:A(\epsilon,\eta)\rightarrow A(\epsilon)$ defined respectively by $p_{\epsilon}(a+\epsilon b+\eta c +\epsilon \eta d)=a+ \eta c$ and $p_{\eta}(a+\epsilon b+\eta c +\epsilon \eta d)=a+ \epsilon b$. Each of these morphisms induces a morphism of group via the functor $\GC$. Then by functoriality we deduce that $\Sigma_{A(\epsilon,\eta)}(Id + \epsilon  M)=Id + \epsilon \Phi_{A}(M)$ and 
$\Sigma_{A(\epsilon,\eta)}(Id + \eta  M)=Id + \eta \Phi_{A}(M)$ for $M \in \GC(A)$.\\
Finally from the following equalities 
\begin{eqnarray}
\Sigma_{A(\epsilon, \eta)}(Id + \epsilon \eta [M,N]_{A}) & = & \Sigma_{A(\epsilon, \eta)}((Id + \epsilon M)(Id + \eta N)(Id - \epsilon M)(Id - \eta N))\nonumber\\
& = & (Id + \epsilon \Phi_{A}(M))(Id + \eta \Phi_{A}(N))(Id - \epsilon \Phi_{A}(M))\nonumber\\
& &(Id - \eta \Phi_{A}(N))\nonumber \\
& = & Id + \epsilon \eta [\Phi_{A}(M),\Phi_{A}(N)]_{A}
\end{eqnarray}
we deduce that $\Phi_{A}$ is a morphism of Lie algebra.
$\qquad\qquad\qquad\qquad\qquad\;\;\blacksquare$
 

\vskip1pc
\noindent We end this section with the following theorem which shows that the Lie superalgebra of the real form $\GC^{\Sigma}$ is the same thing as the real form of $Lie(\GC)$, i.e. $Lie(\GC^{\Sigma})=(Lie(\GC))^{\Phi}$.


\begin{thm}
Let $\Sigma$ be the real structure of the algebraic supergroup $\GC$ and $\Phi$ its corresponding real structure on the Lie superalgebra $Lie(\GC)$. The real form $(Lie(\GC))^{\Phi}$, defined by $(Lie(\GC))^{\Phi}(A)=\{ x \in Lie(\GC), \, \Phi_{A}(x)=x \}$, is then the Lie superalgebra of the real form $\GC^{\Sigma}$.
\end{thm}


\vskip1pc
\noindent
\textbf{Proof:}
By definition, $Lie(\GC^{\Sigma})(A)=\{X \in \GC^{\Sigma}(A(\epsilon)),\; \GC(p)(X)=Id,\; \Sigma_{A}(X)$ $=X \}$. The elements $X$ of $\GC^{\Sigma}(A)$ which fullfil the condition $\GC(p)(X)=Id$ are the even supermatrices of the form $Id+\epsilon M$ with $M \in gl_{m\vert n}(A)$. Therefore,  $Lie(\GC^{\Sigma})(A)=\{M \in gl_{m \vert n}(A), \; \Sigma_{A(\epsilon)}(Id+ \epsilon M)=Id+ \epsilon M \}=\{M \in  M_{m \vert n}(A),\; \Phi_{A}(M)=M \}$ which proves the theorem. 
$\qquad\qquad\qquad\qquad\qquad\blacksquare$



\section{Lifting the Serganova automorphisms to the algebraic supergroups $SL_{m\vert n}$ and $OSp_{m\vert 2n}$}

Before defining the algebraic supergroups $SL_{m\vert n}$ and $OSp_{m\vert 2n}$ we introduce some notations. Let $J_{m,n}$ be the following supermatrix $diag(1_{m},J_{n})$ where
\begin{eqnarray}
diag(M,N)=\left(
\begin{array}{cc}
M & 0 \\
0 & N 
\end{array}
\right)
,
J_{n}=\left(
\begin{array}{cc}
0 & 1_{n} \\
-1_{n} & 0 
\end{array}
\right),
\end{eqnarray}
with $1_{n}$ is the unit square matrix of order $n$.
The supertranspose and the $\Pi$-transpose of an even supermatrix is defined by
\begin{eqnarray}
\left(
\begin{array}{cc}
A & B \\
C & D 
\end{array}
\right)^{st}=\left(
\begin{array}{cc}
A^{t} & -C^{t} \\
B^{t} & D^{t} 
\end{array}
\right),\;\;
\Pi\left(
\begin{array}{cc}
A & B \\
C & D 
\end{array}
\right)=\left(
\begin{array}{cc}
D & C \\
B & A 
\end{array}
\right)
\end{eqnarray}
where $M^{t}$ means the usual transpose of the matrix $M$.
The supertrace of a supermatrix is given by the formula:
\begin{eqnarray}
str\left(
\begin{array}{cc}
A & B \\
C & D 
\end{array}
\right)=trA-trD
\end{eqnarray}
where $trA$ is the usual trace of the matrix $A$.
\vskip1pc
\noindent 
The functor $SL_{m\vert n}$ is defined by:
\begin{eqnarray}
SL_{m\vert n}: \mathcal{A} &  \longrightarrow & \mathcal{S} \\
A  & \longmapsto & SL_{m\vert n}(A)=\{ M \in GL_{m \vert n}(A),\; sdet(M)=1 \}\nonumber
\end{eqnarray}
and the functor $OSp_{m\vert 2n}$ is:
\begin{eqnarray}
OSp_{m\vert 2n}: \mathcal{A} &  \longrightarrow & \mathcal{S} \\
A  & \longmapsto & OSp_{m\vert 2n}(A)=\{ M \in GL_{m \vert 2n}(A),\; M^{st}J_{m,n}M=J_{m,n} \}\nonumber.
\end{eqnarray} 
Their corresponding Lie superalgebras are the functors $sl_{m\vert n}$ and $osp_{m\vert 2n}$ which associate respectively to a superalgebra $A$ the following sets $sl_{m\vert n}(A)=\{X\in gl_{m\vert n}(A), \; str(X)=0\}$ and
$osp_{m\vert 2n}(A)=\{X\in gl_{m\vert 2n}(A),\; X^{st}J_{m,n}+J_{m,n}X=0 \}$.
\vskip1pc
\noindent 
To describe the automorphisms of Serganova for these Lie superalgebras we need to introduce the following conventions:
\begin{eqnarray}
\delta_{\lambda}\left(
\begin{array}{cc}
A & B \\
C & D 
\end{array}
\right)=\left(
\begin{array}{cc}
A & \lambda B \\
\lambda^{-1}C & D 
\end{array}
\right),\;\; 
{I_{n}}^{l}=diag(1_{l},-1_{n-l}),
\end{eqnarray}
\begin{eqnarray}
Ad(M)X=MXM^{-1}, c(N)=\bar{N}
\end{eqnarray}
where $M,X$ are supermatrices of the same order while $N$ is a supermatrix filled with complex numbers and the bar on $N$ indicate that all the entries of $N$ are conjugated. 
\vskip1pc
\noindent
Now, we write the automorphisms of Serganova which give rise to the standard real structure for $sl_{m\vert n}$ (see table 3 in \cite{Serganova}). They are of four types:
\begin{eqnarray}
\sigma_{1}(M)& = &-st \circ Ad(diag({I_{m}}^{p},{I_{m}}^{q})) \circ c \circ \delta_{i}(M)\\
\sigma_{2}(M)& = &Ad(diag(J_{m},J_{n}))\circ c(M) \;\;\;(m,n\vert 2)\\ 
\sigma_{3}(M)& = &\Pi \circ c(M) \;\;\;(m=n) \\
\sigma_{4}(M)& = &-st \circ \Pi \circ c(M) \;\;\;(m=n\vert 2) 
\end{eqnarray}
with $M \in sl_{m\vert n}(\mathbb{C})$. \\
The automorphisms of Serganova which come from graded real structure of $sl_{m \vert n}$ are (see table 6 in \cite{Serganova}):
\begin{eqnarray}
\omega_{1}(M)& = & c \circ Ad(diag(1_{m},J_{n})) \circ (M) \;\;\;(n\vert 2)\\
\omega_{2}(M)& = & -st \circ c \circ Ad(diag({I_{m}}^{p},{I_{n}}^{q}))(M) \\ 
\omega_{3}(M)& = & c \circ \Pi \circ \delta_{i}(M) \;\;\;(m=n) 
\end{eqnarray}
with $M \in sl_{m\vert n}(\mathbb{C})$. \\
In the functorial language, the standard real structures $\bar{\sigma}_{l}$ ($l=1,2,3,4$) associated to the $\sigma_{l}$ are given by the same formula except that the $c$ (the complex conjugation) is replaced by the standard conjugation $a \rightarrow \bar{a}$ of the superalgebra $A$. For the graded real structure $\tilde{\omega}_{k}$ ($k=1,2,3$), the complex conjugation is replaced by the graded conjugation $a \rightarrow \tilde{a}$.
\vskip1pc
\noindent
It turns out that the lifts of the $\tilde{\phi}$ are the following \underline{standard} real structures of the supergroup $SL_{m\vert n}$
\begin{eqnarray}
\Sigma_{1}(X)& = &(-\bar{\sigma}_{1}(X))^{-1} \\
\Sigma_{2}(X)& = &\bar{\sigma}_{2}(X) \;\;\;(m,n\vert 2)\\ 
\Sigma_{3}(X)& = &\bar{\sigma}_{3}(X) \;\;\;(m=n) \\
\Sigma_{4}(X)& = &(-\bar{\sigma}_{4}(X))^{-1} \;\;\;(m=n\vert 2) 
\end{eqnarray}
for $X \in SL_{m \vert n}(A)$.\\
The lifts of the \underline{graded} real structure of the supergroup $SL_{m \vert n}$ are, in turn,
\begin{eqnarray}
\Omega_{1}(X)& = & \tilde{\omega}_{1}(X) \;\;\;(n\vert 2)\\
\Omega_{2}(X)& = & (-\tilde{\omega}_{2}(X))^{-1} \\ 
\Omega_{3}(X)& = & \tilde{\omega}_{3}(X) \;\;\;(m=n) 
\end{eqnarray}
for $X \in SL_{m \vert n}(A)$.


\begin{rem}
It is easy to prove that $\Sigma_{k}$ ($\Omega_{l}$) are indeed the standard (graded) real forms of $SL_{m \vert n}$ in the sense of definition 3.2. Moreover, it is not more difficult to find that these real structures fullfill the following equalities $\Sigma_{l}(Id+ \epsilon M)=Id + \epsilon \bar{\sigma}_{l}(M)$, $\Omega_{l}(Id+ \epsilon M)=Id + \epsilon \tilde{\omega}_{l}(M)$, which means that they are clearly the lifts of $\bar{\sigma}_{k}$, $\tilde{\omega}_{l}$.
\end{rem}


\noindent
Now we turn to the Serganova automorphism of $OSp_{m\vert 2n}$. The Lie algebra automorphisms which come from \underline{standard} real structures are (see table 3  in \cite{Serganova})
\begin{eqnarray}
\xi_{1}(M)& = & Ad(diag({I_{m}}^{p},1_{2n}))\circ c(M) \\
\xi_{2}(M)& = & Ad(diag(J_{m},{I_{n}}^{p},{I_{n}}^{p}))(M)  \;\;\;(m\vert 2)
\end{eqnarray}
for $M\in osp_{m\vert 2n}(\mathbb{C})$.\\ 
The Lie algebra automorphisms which give rise to the \underline{graded} real structure are (see table 6  in \cite{Serganova})
\begin{eqnarray}
\psi_{1}(M)& = & c \circ Ad(diag({I_{m}}^{p},d({I_{n}}^{q},{I_{n}}^{q}))\circ J_{2n})\circ (M) \\
\psi_{2}(M)& = & c \circ Ad(diag(J_{m},1_{2n}))\circ (M)  \;\;\;(m\vert 2)
\end{eqnarray}
for $M\in osp_{m\vert 2n}(\mathbb{C})$.\\
In order to switch to the functorial language, we again replace the complex conjugation in the automorphisms $\xi_{1,2}$ ($\psi_{1,2}$) by the standard (graded) conjugation and we denote the corresponding standard (graded) real structure  by $\bar{\xi}_{1,2}$ ($\tilde{\psi}_{1,2}$).\\ 
Their lifts to the supergroup $OSp_{m\vert n}$ are respectively as follows
\begin{eqnarray}
\Xi_{1}(X)& = & \bar{\xi}_{1}(X) \\
\Xi_{2}(X)& = & \bar{\xi}_{2}(X)  \;\;\;(m\vert 2)\\
\Psi_{1}(X)& = & \tilde{\psi}_{1}(X) \\
\Psi_{2}(X)& = & \tilde{\psi}_{2}(X)  \;\;\;(m\vert 2)
\end{eqnarray}
for $X \in OSp_{m,2n}(A)$. 


\begin{rem}
it is easy to prove that $\Xi_{k}$ ($\Psi_{l}$) are indeed the standard (graded) real forms of $OSp_{m \vert 2n}$ in the sense of definition 3.2. Moreover, it is not more difficult to find that these real structures fullfill the following equalities $\Xi_{l}(Id+ \epsilon M)=Id + \epsilon \bar{\xi}_{l}(M)$, $\Psi_{l}(Id+ \epsilon M)=Id + \epsilon \tilde{\psi}_{l}(M)$, which means that they are clearly the lifts of $\bar{\xi}_{k}$, $\tilde{\psi}_{l}$.
\end{rem}


 \end{document}